\documentclass[12pt]{amsart}
\usepackage{fullpage}
\DeclareFontEncoding{OT2}{}{} % to enable usage of cyrillic fonts

\newcommand{\bad}{{\operatorname{bad}}}
\newcommand{\calP}{{\mathcal P}}

\newcommand{\Q}{{\mathbb Q}}
\newcommand{\Qbar}{{\overline{\Q}}}

\newcommand{\Z}{{\mathbb Z}}
\newcommand{\R}{{\mathbb R}}

\newcommand{\PP}{{\mathbb P}}
\newcommand{\OO}{{\mathcal O}}

\newcommand{\F}{{\mathbb F}}

\newcommand{\Aut}{\operatorname{Aut}}
\newcommand{\Gal}{\operatorname{Gal}}

\newcommand{\Spec}{\operatorname{Spec}}

\newcommand{\GL}{{\operatorname{GL}}}

\newcommand{\isom}{\simeq}

\newtheorem{theorem}{Theorem}[section]
\newtheorem{lemma}[theorem]{Lemma}
\newtheorem{corollary}[theorem]{Corollary}
\newtheorem{proposition}[theorem]{Proposition}

\theoremstyle{definition}
\newtheorem{definition}[theorem]{Definition}

\theoremstyle{remark}
\newtheorem{remark}[theorem]{Remark}

\begin{document}

\title[Large subrings of $\Q$]{Hilbert's Tenth Problem and Mazur's Conjecture for large subrings of $\Q$}
\author{Bjorn Poonen}
\thanks{This research was supported by NSF grant DMS-0301280,
	and a Packard Fellowship.  The paper will appear 
	in {\em J.\ Amer.\ Math.\ Soc.}}
\address{Department of Mathematics, University of California, Berkeley, CA 94720-3840, USA}
\email{poonen@math.berkeley.edu}
\date{June 14, 2003}

\begin{abstract}
We give the first examples of infinite sets of primes $S$
such that Hilbert's Tenth Problem over $\Z[S^{-1}]$ 
has a negative answer.
In fact, we can take $S$ to be a density~1 set of primes.
We show also that for some such $S$ 
there is a punctured elliptic curve $E'$ over $\Z[S^{-1}]$
such that the topological closure of $E'(\Z[S^{-1}])$ in $E'(\R)$
has infinitely many connected components.
\end{abstract}

\maketitle

%****************************************************************************
\section{Introduction}
\label{introduction}

Hilbert's Tenth Problem, in modern terms, 
was to find an algorithm (Turing machine) to decide,
given a polynomial equation $f(x_1,\dots,x_n)=0$
with coefficients in $\Z$, 
whether it has a solution with $x_1,\dots,x_n \in \Z$.
Y.~Matijasevi{\v{c}}~\cite{matijasevic1970}, 
building on earlier work of M.~Davis, H.~Putnam, 
and J.~Robinson~\cite{davis-putnam-robinson1961},
showed that no such algorithm exists.
If one replaces $\Z$ in both places by a different 
commutative ring $R$ (let us assume its elements can be and have been encoded 
for input into a Turing machine),
one obtains a different question,
called {\em Hilbert's Tenth Problem over $R$},
whose answer depends on $R$.
These problems are discussed in detail in~\cite{hilbertstenthproblem}.

In particular, the answer for $R=\Q$ is unknown.
Hilbert's Tenth Problem over $\Q$ is equivalent 
to the general problem of deciding whether a variety over $\Q$ 
has a rational point.
One approach to proving that Hilbert's Tenth Problem over $\Q$ 
has a negative answer would be to deduce this from
Matijasevi{\v{c}}'s theorem for $\Z$,
by showing that $\Z$ is diophantine over $\Q$ in the 
following sense:

\begin{definition}
Let $R$ is a ring, and $A \subseteq R^m$.
Then $A$ is {\em diophantine over $R$} if and only if there exists
a polynomial $f$ in $m+n$ variables with coefficients in $R$ such that
	$$A = \{\, a \in R^m \mid \exists x \in R^n 
			\text{ such that } f(a,x)=0\,\}.$$
\end{definition}

On the other hand, Mazur conjectures that if $X$ is a
variety over $\Q$, then the topological closure of $X(\Q)$
in $X(\R)$ has only finitely many 
components~\cite{mazur1992},\cite{mazur1995}.
This would imply that $\Z$ is not diophantine over $\Q$.
More generally, Cornelissen and Zahidi~\cite{cornelissen-zahidi2000} 
have shown that Mazur's Conjecture
implies that there is no diophantine model of $\Z$ over $\Q$.

\begin{definition}
A {\em diophantine model} of $\Z$ over $\Q$
is a set $A \subseteq \Q^n$ that is diophantine over $\Q$
with a bijection $\Z \to A$
under which the graphs of addition and multiplication on $\Z$
correspond to subsets of $A^3 \subseteq \Q^{3n}$
that are diophantine over $\Q$.
\end{definition}

This is important,
because the existence of such a diophantine model,
together with Matijasevi{\v{c}}'s Theorem,
would imply a negative answer for Hilbert's Tenth Problem over $\Q$.

\medskip

This paper studies Hilbert's Tenth Problem over rings
between $\Z$ and $\Q$.
Such rings are in bijection with subsets of the set $\calP$
of prime numbers.
Namely, given $S \subseteq \calP$, one has the ring $\Z[S^{-1}]$,
and conversely, given a subring $R$ between $\Z$ and $\Q$,
one has $R=\Z[S^{-1}]$ where $S=\calP \cap R^\times$.

Using quadratic forms as in J.~Robinson's work, 
one can show that for any prime $p$ the ring $\Z_{(p)}$
of rational numbers with denominators prime to $p$
is diophantine over $\Q$~\cite[Proposition~3.1]{kim-roush1992b}.
{}A short argument using this shows that for finite $S$,
Hilbert's Tenth Problem over $\Z[S^{-1}]$ has a negative answer.

In this paper, we give the first examples of \emph{infinite} subsets $S$
of $\calP$ for which Hilbert's Tenth Problem over $\Z[S^{-1}]$
has a negative answer.
In fact, we show that there exist such $S$ of natural density $1$,
so in one sense, we are approaching a negative answer for $\Q$.
(See Section~\ref{densities} for the definition of natural density.)
Previously, Shlapentokh proved that if $K$ is a totally real number field
or a totally complex degree-$2$ extension of a totally real number field,
then there exists a set of places $S$ of $K$ 
of Dirichlet density arbitrarily close to $1 - [K:\Q]^{-1}$
such that if $\OO_{K,S}$ is the subring of elements of $K$
that are integral at all places outside $S$,
then Hilbert's Tenth Problem over $\OO_{K,S}$ has a negative 
answer~\cite{shlapentokh1997},\cite{shlapentokh2000highdensity},\cite{shlapentokh2002}.
But for $K=\Q$, 
this gives nothing beyond Matijasevi{\v{c}}'s Theorem.

More generally, we prove the following:

\begin{theorem}
\label{main}
There exist disjoint recursive sets of primes $T_1$ and $T_2$,
both of natural density $0$, such that 
for any set $S$ of primes containing $T_1$ and disjoint from $T_2$,
the following hold:
\begin{enumerate}
\item\label{firstpart}
There exists an affine curve $E'$ over $\Z[S^{-1}]$
such that the topological closure of $E'(\Z[S^{-1}])$ in $E'(\R)$
is an infinite discrete set.
\item\label{secondpart}
The set of positive integers with addition and multiplication
admits a diophantine model over $\Z[S^{-1}]$.
\item\label{thirdpart}
Hilbert's Tenth Problem over $\Z[S^{-1}]$ has a negative answer.
\end{enumerate}
\end{theorem}

\begin{remark} \hfill
\begin{enumerate}
\item[(i)]
Arguably~\eqref{thirdpart} is the most important of the 
three parts.  
We have listed the parts in the order they will be proved.
\item[(ii)] 
A subset $T \subseteq \Z$ is {\em recursive} if and only if
there exists an algorithm (Turing machine) 
that takes as input an integer $t$ 
and outputs YES or NO according to whether $t \in T$.
\item[(iii)] 
We use natural density instead of Dirichlet density
in order to have a slightly stronger statement.
See~\cite[VI.4.5]{serrecourse} for the definition
of Dirichlet density and its relation to natural density.
\end{enumerate}
\end{remark}

Previously, Shlapentokh~\cite{shlapentokh2003} used norm equations 
to prove that there
exist sets $S \subseteq \calP$ of Dirichlet density arbitrarily close to $1$
for which there exists an affine variety $X$ over $\Z[S^{-1}]$
such that the closure of $X(\Z[S^{-1}])$ in $X(\R)$
has infinitely many connected components.
(She also proved an analogous result for localizations of the
ring of integers of totally real number fields and totally complex
degree-$2$ extensions of totally real number fields.
For number fields with exactly one conjugate pair of nonreal embeddings,
she obtained an analogous result, but with density only $1/2$.)

Question~4.1 of~\cite{shlapentokh2003}
asked whether over $\Q$ one could do the same for some $S \subseteq \calP$ 
of Dirichlet density exactly $1$.
Part~(1) of our Theorem~\ref{main} 
gives an affirmative answer (take $S=\calP - T_2$).
In fact, it was the attempt to answer Shlapentokh's Question~4.1 that
inspired this paper, so the author thanks her
for asking the right question.

The rest of this paper is devoted to proving Theorem~\ref{main}.
The strategy is to take an elliptic curve $E$ over $\Q$ such
that $E(\Q)$ is generated
by one point $P$ of infinite order,
and to construct $T_1$ (resp.\ $T_2$) so that certain prime multiples
(resp.\ at most finitely many other integer multiples)
of $P$ have coordinates in $\Z[S^{-1}]$.
Using Vinogradov's result on the equidistribution 
of the prime multiples of an irrational number modulo~1,
we can prescribe the approximate locations of the prime multiples 
of $P$ in $E(\R)$.
If we prescribe them so that their $y$-coordinates approximate
the set of positive integers sufficiently well, 
then approximate addition and approximate squaring
on the set $A$ of these $y$-coordinates 
make $A$ into a diophantine model of the positive integers.

\medskip

Shlapentokh and the author plan eventually to write a joint paper
generalizing Theorem~\ref{main} to other number fields,
and to places other than the real place.

%****************************************************************************
\section{Elliptic curve setup}
\label{setup}

Let $E$ be an elliptic curve over $\Q$ of rank~1.
To simplify the arguments, we will assume moreover that $E(\Q) \isom \Z$,
that $E(\R)$ is connected, 
and that $E$ does not have complex multiplication.
For example, these conditions hold for the smooth projective model
of $y^2=x^3+x+1$.
Let $P$ be a generator of $E(\Q)$.
Fix a Weierstrass equation $y^2=x^3+ax+b$ for $E$, where $a,b \in \Z$.

Let $E'=\Spec \Z[S_\bad^{-1}][x,y]/(y^2 - (x^3+ax+b))$,
where $S_\bad$ is a finite set of primes such that 
$E'$ is smooth over $\Z[S_\bad^{-1}]$.
In particular, $2 \in S_\bad$.
Enlarge $S_\bad$ if necessary so that $P \in E'(\Z[S_\bad^{-1}])$.

%****************************************************************************
\section{Denominators of $x$-coordinates}
\label{S:denominators}

For nonzero $n \in \Z$, let $d_n \in \Z_{>0}$
be the prime-to-$S_\bad$ part of
the denominator of $x(nP)$;
that is, $d_n$ is the product one obtains
if one takes the prime factorization of the denominator
of $x(nP)$ and omits the powers of primes in $S_\bad$.
Define $d_0=0$.
The notation $m \mid n$ means $n \in m\Z$.

\begin{lemma} \hfill
\label{denominators}
\begin{enumerate}
\item[(a)]
For any $r \in \Z$, the set $\{\, n \in \Z : r \mid d_n\,\}$
is a subgroup of $\Z$.
\item[(b)]
There exists $c \in \R_{>0}$ such that 
$\log d_n = (c-o(1)) n^2$ as $n \to \infty$ (cf.~\cite[Lemma~8]{silverman1988}).
\end{enumerate}
\end{lemma}

\begin{proof}
(a) 
We may reduce to the case where $r=p^e$ for some
prime $p \not\in S_\bad$ and $e \in \Z_{>0}$.
Thus it suffices to show that 
the set $E_e:=\{\, Q \in E(\Q_p): v_p(x(Q)) \le -e \,\} \cup \{O\}$
is a subgroup of $E(\Q_p)$,
where $v_p \colon \Q_p^* \to \Z$ is the $p$-adic valuation.
The set $E_1$ is the kernel of the reduction map $E(\Q_p) \to E(\F_p)$
(extend $E$ to a smooth proper curve over $\Z_p$ to make sense of this).
Since $p>2$, the formal logarithm 
$\lambda \colon E_1 \to p\Z_p$ is an isomorphism \cite[IV.6.4]{silvermanAEC}.
By \cite[IV.5.5,~IV.6.3]{silvermanAEC},
$v_p(\lambda(Q))=v_p(z(Q))$ for all $Q \in E_1$,
where $z=-x/y$ is the standard parameter for the formal group.
By \cite[pp.~113--114]{silvermanAEC}, 
$x = -2 z^{-2} + \dots \in \Z_p((z))$,
so $v_p(x(Q))=-2 v_p(z(Q)) = -2 v_p(\lambda(Q))$.
Thus $G_e$ corresponds under $\lambda$ to $p^{\lceil e/2 \rceil} \Z_p$,
and is hence a subgroup.

(b) The number $\log d_n$ is the logarithmic height $h(nP)$,
except that in the sum defining the height,
the terms corresponding to places in $S_\bad \cup \{\infty\}$
have been omitted.
A standard diophantine approximation result
(see Section~7.4 of~\cite{serremordellweil})
implies that each such term 
contributes at most a fraction $o(1)$ of the height,
as $n \to \infty$.
If $\hat{h}$ is the canonical height, then
$h(nP)=\hat{h}(nP) + O(1) = \hat{h}(P) n^2 + O(1)$.
Take $c=\hat{h}(P)$, which is positive, since $P$ is not torsion.
\end{proof}

\begin{remark}
The bottom of p.~306 in~\cite{ayad1992}
relates the denominators of $x(nP)$ in lowest terms 
to the sequence of values of division polynomials evaluated at $P$.
The study of divisibility properties of the latter sequence
is very old: results were claimed in the $19^{\operatorname{th}}$ century 
by Lucas (but apparently not published),
and proofs were given in~\cite{ward1948}.
\end{remark}

For $n \in \Z$, let $S_n$ be the set of prime factors of $d_n$.
If $m,n \in \Z$, then $(m,n)$ denotes their greatest common divisor.

\begin{corollary}
\label{cor:gcd}
If $m,n \in \Z$, then $S_{(m,n)} = S_m \cap S_n$.
In particular, if $(m,n)=1$, then $S_m$ and $S_n$ are disjoint.
\end{corollary}

\begin{proof}
Lemma~\ref{denominators}(a) implies the first statement.
Since $P \in E'(\Z[S_\bad^{-1}])$, we have $S_1=\emptyset$,
and the second statement follows.
\end{proof}

\begin{lemma}
\label{setdifference}
If $\ell$ and $m$ are primes, and $\max\{\ell,m\}$ is sufficiently large,
then $S_{\ell m} - (S_\ell \cup S_m)$ is nonempty.
\end{lemma}

\begin{proof}
If $p \mid d_m$, or equivalently $v_p(x(mP)) < 0$,
then using the formal logarithm $\lambda$ 
as in the proof of Lemma~\ref{denominators}(a)
we obtain
\[
	v_p(d_{\ell m}) = -v_p(x(\ell m P)) = 2 v_p(\lambda(\ell m P))
	= 2 v_p(\ell \lambda(mP)) = v_p(\ell^2 d_m).
\]
If $S_{\ell m} - (S_\ell \cup S_m)$ were empty,
then for each $p\mid d_{\ell m}$
we could apply either this result or the analogue with 
$\ell$ and $m$ interchanged,
and hence deduce $d_{\ell m} \mid \ell^2 m^2 d_\ell d_m$.
This contradicts Lemma~\ref{denominators}(b) if $\max\{\ell,m\}$
is sufficiently large.
\end{proof}

\begin{remark}
Our Lemma~\ref{setdifference} is a special case of
Lemma~9 of~\cite{silverman1988} 
(except for the minor differences that \cite{silverman1988}
requires $E$ to be in minimal Weierstrass form
and considers the full denominator instead of its
prime-to-$S_\bad$ part). 
The method of proof is the same.
These results may viewed as elliptic analogues 
of Zsigmondy's Theorem: 
see~\cite{everest2002}.
\end{remark}

%****************************************************************************
\section{Definition of $T_1$ and $T_2$}

For each prime number $\ell$, let $a_\ell$ be the smallest 
$a \in \Z_{>0}$ such that $d_{\ell^a} > 1$.
By Lemma~\ref{denominators}(b), $a_\ell$ exists,
and $a_\ell=1$ for all $\ell$ outside a finite set $L$ of primes.
Baker's method~\cite[Chapter~8]{serremordellweil} 
lets us compute the finite set $E'(\Z[S_\bad^{-1}])$,
so the set $L$ and the values $a_\ell$ for $\ell \in L$ are computable.

Let $p_\ell = \max S_{\ell^a}$ where $a=a_\ell$.
For primes $\ell$ and $m$ (possibly equal),
Lemma~\ref{setdifference} lets us define 
$p_{\ell m} = \max \left(S_{\ell m} - (S_\ell \cup S_m) \right)$
when $\max\{\ell,m\}$ is sufficiently large.
Let $\ell_1<\ell_2<\dots$ be a sequence of primes outside $L$.
(The $\ell_i$ will be constructed in Section~\ref{primeconstruction} 
with certain properties, 
but for now these properties are not relevant.)
The plan will be to force $\ell_i P \in E'(\Z[S^{-1}])$ for all $i$,
by requiring each $S_{\ell_i}$ to be contained in $S$.
On the other hand, we must require other primes to lie outside $S$
to make sure that not too many other multiples of $P$
end up in $E'(\Z[S^{-1}])$.

Let $T_1 = S_\bad \cup \bigcup_{i \ge 1} S_{\ell_i}$.
Let $T_2^a$ be the set of $p_\ell$
for all primes $\ell \not\in \{\ell_1,\ell_2,\dots\}$.
If $\ell_1$ is sufficiently large, we may define 
$T_2^b = \{\, p_{\ell_i \ell_j} : 1 \le j \le i \,\}$
and $T_2^c = \{\, p_{\ell \ell_i} : \ell \in L, i \ge 1 \,\}$.
Finally, let $T_2 = T_2^a \cup T_2^b \cup T_2^c$.

%****************************************************************************
\section{Properties of $T_1$ and $T_2$}

\begin{proposition}
The sets $T_1$ and $T_2$ are disjoint.
\end{proposition}

\begin{proof}
By definition of $p_\ell$ and $p_{\ell m}$,
$S_\bad \cap T_2 = \emptyset$.
If $\ell \not= \ell_i$, then $(\ell^{a_\ell},\ell_i)=1$,
so Corollary~\ref{cor:gcd} implies $p_\ell \not\in S_{\ell_i}$.
Thus $T_1 \cap T_2^a =\emptyset$.
If $i \not\in \{j,k\}$,
then Corollary~\ref{cor:gcd} implies $p_{\ell_j \ell_k} \not\in S_{\ell_i}$,
while if $i \in \{j,k\}$ then $p_{\ell_j \ell_k} \not\in S_{\ell_i}$
by definition of $p_{\ell_j \ell_k}$.
Thus $T_1 \cap T_2^b =\emptyset$.
If $i \not=j$,
then Corollary~\ref{cor:gcd} implies $p_{\ell \ell_j} \not\in S_{\ell_i}$,
while if $i=j$, then $p_{\ell \ell_i} \not\in S_{\ell_i}$
by definition of $p_{\ell_j \ell_k}$.
Thus $T_1 \cap T_2^c =\emptyset$.
\end{proof}

\begin{proposition}
\label{integerpoints}
If $S$ contains $T_1$ and is disjoint from $T_2$,
then $E'(\Z[S^{-1}])$ is the union of $\{\,\pm \ell_i P: i \ge 1 \,\}$
and some subset of the finite set
$\left\{\, rP : r \mid \prod_{\ell \in L} \ell^{a_\ell-1} \,\right\}$.
\end{proposition}

\begin{proof}
Because the equation of $E$ relates the $x$- and $y$-coordinates,
a point $nP$ belongs to $E'(\Z[S^{-1}])$
if and only if $S_n \subseteq S$.
In particular, $S_{\ell_i} \subseteq T_1 \subseteq S$,
so $\pm \ell_i P \in E'(\Z[S^{-1}])$.

Any point outside 
	$$\{\,\pm \ell_i P: i \ge 1 \,\} \cup
	\left\{\, rP : r \mid \textstyle 
		\prod_{\ell \in L} \ell^{a_\ell-1} \,\right\}$$
is $nP$ for some $n$ divisible by one of the following:
\begin{itemize} 
\item $\ell^{a_\ell}$ for some $\ell$ 
	not in the sequence $\ell_1,\ell_2,\dots$,
\item $\ell_i \ell_j$ for some $1 \le j \le i$, or
\item $\ell \ell_i$ for some $\ell \in L$ and $i \ge 1$.
\end{itemize}
Lemma~\ref{denominators}(a) implies then that $S_n$
contains a prime of $T_2^a$, $T_2^b$, or $T_2^c$,
respectively, so $S_n \not\subseteq S$.
\end{proof}

%****************************************************************************
\section{Natural density}
\label{densities}

The {\em natural density} of a subset $T \subseteq \calP$ is defined as 
	$$\lim_{X \to \infty} \frac{\# \{p \in T: p \le X\}}
	{\# \{p \in \calP: p \le X\}},$$
if the limit exists.
One defines {\em upper natural density} similarly, using $\limsup$ instead of $\lim$.

\begin{lemma}
\label{vinogradov}
If $\alpha \in \R-\Q$, then 
$\{\, \ell \alpha \bmod 1 : \ell \text{ is prime}\,\}$
is equidistributed in $[0,1]$.
That is, for any interval $I \subseteq [0,1]$,
the set of primes $\ell$ for which $(\ell \alpha \bmod 1)$
belongs to $I$ has natural density equal to the length of $I$.
\end{lemma}

\begin{proof}
See p.~180 of~\cite{vinogradov1954}.
\end{proof}

Let $y(\ell P) \in \Q$ denote the $y$-coordinate of $\ell P \in E(\Q)$.

\begin{corollary}
\label{ycoordinates}
If $I \subseteq \R$ is an interval with nonempty interior,
then the set of primes $\ell$ for which $y(\ell P) \in I$
has positive natural density.
\end{corollary}

\begin{proof}
Since $E(\R)$ is a connected compact $1$-dimensional Lie group over $\R$, 
we can choose an isomorphism $E(\R) \to \R/\Z$ as topological groups.
Since $P$ is of infinite order, its image in $\R/\Z$
is represented by an irrational number.
The subset of $E(\R)$ having $y$-coordinate in $I$ 
corresponds to a nontrivial interval in $\R/\Z$.
Now apply Lemma~\ref{vinogradov}.
\end{proof}

%****************************************************************************
\section{Construction of the $\ell_i$}
\label{primeconstruction}

For prime $\ell$, define
	 $$\mu_\ell = \sup_{X \in \Z_{\ge 2}} 
		\frac{\# \{p \in S_\ell : p \le X\}}
		{\# \{p \in \calP: p \le X\}}.$$
The supremum is attained for some $X \le \max S_\ell$,
so $\mu_\ell$ is computable for each $\ell$.

\begin{lemma}
\label{mulemma}
For any $\epsilon>0$,
the natural density of $\{\, \ell : \mu_\ell > \epsilon \,\}$ is $0$.
\end{lemma}

\begin{proof}
For $X \in \R$, let $\pi(X):=\# \{p \in \calP: p \le X\}$.
If $\ell$ is a prime and $\mu_\ell>\epsilon$,
then we can choose $X_\ell \in \Z_{\ge 2}$ such that
	$$\frac{\# \{p \in S_\ell : p \le X_\ell\}}
	{\pi(X_\ell)} > \epsilon.$$

For $M \in \Z_{\ge 2}$, let $U_M$ be the set of primes $\ell$
such that $\mu_\ell >\epsilon$ and $X_\ell \in [M,2M)$.
If $\ell \in U_M$, then
	$$\# \{p \in S_\ell : p \le 2M\} 
		\quad \ge \quad
	\# \{p \in S_\ell : p \le X_\ell\} 
		\quad > \quad 
	\epsilon \pi(X_\ell) 
		\quad \ge \quad
	\epsilon \pi(M).$$
But the $S_\ell$ are disjoint by Corollary~\ref{cor:gcd}, so
	$$\pi(2M) 
		\quad \ge \quad 
	\sum_{\ell \in U_M} \# \{p \in S_\ell : p \le 2M\}
		\quad \ge \quad 
	\epsilon \pi(M) \#U_M.$$
Thus by the Prime Number Theorem, $\#U_M=O(1)$ as $M \to \infty$.
If $2^{k-1} \le N < 2^k$, then
	$$\#\{\, \ell : \mu_\ell > \epsilon \text{ and } X_\ell \le N \,\} 
	\le \#U_2 + \#U_4 + \#U_8 + \dots + \#U_{2^k} = O(k) = O(\log N)$$
as $N \to \infty$.
If $\mu_\ell>\epsilon$ then by definition of $X_\ell$,

	$$\pi(X_\ell) < \frac{ \# S_\ell}{\epsilon} 
	\le \frac{\log_2 d_\ell}{\epsilon} = O(\ell^2)$$
as $\ell \to\infty$ by Lemma~\ref{denominators}(b),
so $X_\ell = O(\ell^2 \log \ell)$ by the Prime Number Theorem.
Combining the previous two sentences shows that 
\begin{align*}
	\#\{\, \ell\le Y : \mu_\ell > \epsilon \,\} 
	&= \#\{\, \ell\le Y : \mu_\ell > \epsilon \text{ and } 
			X_\ell \le O(Y^2 \log Y) \,\} \\
	&= O(\log O(Y^2 \log Y)),
\end{align*}
which is $o(\pi(Y))$ as $Y \to \infty$.
\end{proof}

Define the $\ell_i$ inductively as follows.
Given $\ell_1$, \dots, $\ell_{i-1}$,
let $\ell_i$ be the smallest prime outside $L$
such that all of the following hold:
\begin{enumerate}
	\item $\ell_i > \ell_j$ for all $j<i$,
	\item $\mu_{\ell_i} \le 2^{-i}$,
	\item $p_{\ell_i \ell_j} > 2^i$ for all $j \le i$,
	\item $p_{\ell \ell_i} > 2^i$ for all $\ell \in L$, and
	\item $|y(\ell_i P) - i| \le 1/(10i)$.
\end{enumerate}

\begin{proposition}
The sequence $\ell_1,\ell_2,\dots$ is well-defined and computable.
\end{proposition}

\begin{proof}
By induction, we need only show that for each $i$,
there exists $\ell_i$ as above.
By Corollary~\ref{ycoordinates}, 
the set of primes satisfying~(5) has positive natural density.
By Lemma~\ref{mulemma}, (2) fails for a set of natural density $0$.
Therefore it will suffice to show that (1), (3), and~(4) are satisfied
by all sufficiently large $\ell_i$.

For fixed $j \le i$, the primes $p_{\ell_i \ell_j}$ for varying
values of $\ell_i$ are distinct by Corollary~\ref{cor:gcd},
so eventually they are greater than $2^i$.
The same holds for $p_{\ell \ell_i}$ for fixed $\ell \in L$.
Thus by taking $\ell_i$ sufficiently large,
we can make all the $p_{\ell_i \ell_j}$ and $p_{\ell \ell_i}$
greater than $2^i$.

Each $\ell_i$ can be computed by searching primes in increasing order
until one is found satisfying the conditions.
\end{proof}

%****************************************************************************
\section{Recursiveness of $T_1$ and $T_2$}

The set $\{\ell_1,\ell_2,\dots\}$
is recursive, since it is a strictly increasing sequence
whose terms can be computed in order.
This is needed for the proofs in this section.

\begin{proposition}
The set $T_1$ is recursive.
\end{proposition}

\begin{proof}
Since $S_\bad$ is finite,
it suffices to give an algorithm for deciding whether 
a prime $p \not\in S_\bad$
belongs to $\bigcup_{i \ge 1} S_{\ell_i}$.
We have $p \in \bigcup_{i \ge 1} S_{\ell_i}$
if and only if $p \mid d_{\ell_i}$ for some $i$,
which holds if and only if the order $n_p$ of the image of $P$ in $E(\F_p)$
divides $\ell_i$ for some $i$.
The order $n_p$ can be computed, and $n_p \not=1$,
since $P \in E'(\Z[S_\bad^{-1}])$.
So we simply check whether $n_p \in \{\ell_1,\ell_2,\dots\}$.
\end{proof}

\begin{lemma}
\label{exactorder}
If $\ell$ is prime, then $\ell \mid \#E(\F_{p_\ell})$.
\end{lemma}

\begin{proof}
By definition of $p_\ell$,
the point $\ell^{a_\ell} P$ reduces to $0$ in $E(\F_{p_\ell})$
but $\ell^{a_\ell-1} P$ does not.
\end{proof}

\begin{proposition}
The set $T_2^a$ is recursive.
\end{proposition}

\begin{proof}
If $p \in T_2^a$,
then $p=p_\ell$ for some $\ell \not\in \{\ell_1,\ell_2,\dots\}$,
and then $\ell \mid \#E(\F_p)$ by Lemma~\ref{exactorder}.
Therefore to test whether a prime $p \not\in S_\bad$ belongs to $T_2^a$,
compute $\#E(\F_p)$ and its prime factors: 
one has $p \in T_2^a$ if and only if there is a prime factor $\ell$
such that $\ell \not\in \{\ell_1,\ell_2,\dots\}$
and $p_\ell=p$.
\end{proof}

\begin{proposition}
The sets $T_2^b$ and $T_2^c$ are recursive.
\end{proposition}

\begin{proof}
By condition~(3) in the definition of $\ell_i$,
if a prime $p$ belongs to $T_2^b$,
it must equal $p_{\ell_i \ell_j}$ for some $1\le j \le i$
with $2^i < p$.
Thus to test whether a prime $p$ belongs to $T_2^b$, 
simply compute $p_{\ell_i \ell_j}$
for $1 \le j \le i < \log_2 p$.

The proof that $T_2^c$ is recursive is similar, using condition~(4).
\end{proof}

Thus $T_1$ and $T_2$ are recursive.

%****************************************************************************
\section{The densities of $T_1$ and $T_2$}

\begin{proposition}
\label{T_1}
The set $T_1$ has natural density $0$.
\end{proposition}

\begin{proof}
For fixed $r \in \Z_{>0}$,
the set $\bigcup_{i>r} S_{\ell_i}$ differs from $T_1$
in only finitely many primes,
so it suffices to show that the former has upper natural density 
tending to $0$ as $r \to \infty$.
By definition of $\mu_{\ell_i}$,
the upper natural density is bounded by 
$\sum_{i>r} \mu_{\ell_i} \le \sum_{i > r} 2^{-i} = 2^{-r}$,
which tends to $0$ as $r \to \infty$.
\end{proof}

\begin{proposition}
\label{T_2bc}
The sets $T_2^b$ and $T_2^c$ have natural density $0$.
\end{proposition}

\begin{proof}
Suppose $2^m \le X < 2^{m+1}$.
By condition~(3) defining $\ell_i$,
the only primes of the form $p_{\ell_i \ell_j}$ 
that might be $\le X$ are those with $1 \le j \le i \le m$.
There are at most $O(m^2)=O((\log X)^2)$ of these, 
which is negligible compared to $\pi(X)$.
Thus $T_2^b$ has natural density $0$.
The proof for $T_2^c$ is similar.
\end{proof}

The rest of this section is devoted to the proof
that $T_2^a$ has natural density $0$.
Recall that $T_2^a$ consists of primes of the form $p_\ell$.
If the sequence of $p_\ell$ grew faster 
than the sequence of primes $\ell$,
then $T_2^a$ would have density $0$.
But Lemma~\ref{exactorder} implies only that
$p_\ell$ is at least about the size of $\ell$.
The strategy for strengthening this bound
will be to show that numbers of the form
$\#E(\F_p)$ are typically divisible by
many primes.
For $n \in \Z_{>0}$, 
let $\omega(n)$ be the number of distinct prime factors of $n$.

\begin{lemma}
\label{E(F_p)}
For any $t \ge 1$, 
the natural density of $\{\, p : \omega(\#E(\F_p)) < t \,\}$ is $0$.
\end{lemma}

\begin{proof}
For a prime $\ell$,
let $E[\ell]$ denote the group of points of order dividing $\ell$ on $E$.
Then $\ell \mid \#E(\F_p)$ if and only if the image of
the Frobenius element at $p$
under $\Gal(\Qbar/\Q) \to \Aut E[\ell]$
has a nonzero fixed vector.
Since $E$ does not have complex multiplication,
the image of $\Gal(\Qbar/\Q)$ in 
$\prod_\ell \Aut E[\ell]$ is open. 
(This follows from~\cite{serre1972}.) 
Thus $\Gal(\Qbar/\Q) \to \prod_{\ell \not\in L'} \Aut E[\ell]$
is surjective for some finite $L' \subseteq \calP$.
A calculation shows that
the fraction of elements of $\Aut E[\ell] \isom \GL_2(\F_\ell)$
having a nonzero fixed vector is 
	$$\frac{\ell^3-2 \ell}{(\ell^2-1)(\ell^2-\ell)} 
	= \frac{1}{\ell} + O\left(\frac{1}{\ell^2}\right).$$
The sum of this over $\ell$ diverges,
so as $C \to \infty$,
the fraction of elements of 
$\prod_{\ell<C, \ell \not\in L'} \Aut E[\ell]$
having fewer than $t$ components with a nonzero fixed vector
tends to $0$.
Applying the Chebotarev Density Theorem 
(see Th\'eor\`eme~1 of~\cite{serre1981} for a version using
natural density) 
and letting $C \to \infty$, we obtain the result.
\end{proof}

\begin{proposition}
\label{T_2a}
The set $T_2^a$ has natural density $0$.
\end{proposition}

\begin{proof}
Because of Lemma~\ref{E(F_p)},
it suffices to show that the upper natural density
of 
\[
	T_2^{a,t}:= \{\, p \in T_2^a : \omega(\#E(\F_p)) \ge t \,\}
\]
tends to $0$ as $t \to \infty$.

Suppose $p = p_\ell \in T_2^{a,t}$.
By Lemma~\ref{exactorder}, $\ell \mid \#E(\F_p)$.
By definition of $T_2^{a,t}$, 
the integer $\#E(\F_p)$ is divisible by at least $t-1$ other primes,
so $2^{t-1} \ell \le \#E(\F_p)$.
There exists a degree-$2$ map $E \to \PP^1$ over $\F_p$,
so $\#E(\F_p) \le 2(p+1) \le 4p$.
Combining the previous two sentences yields 
$\ell \le 2^{3-t} p$.
Since every element of $T_2^{a,t}$ is $p_\ell$ for some $\ell$,
we have
	$$\# \{\, p \in T_2^{a,t} : p \le X \,\} \le \pi(2^{3-t} X) 
	= (2^{3-t} + o(1)) \, \pi(X)$$
as $X \to \infty$.
Thus by definition, the upper natural density of $T_2^{a,t}$
is at most $2^{3-t}$.
This goes to $0$ as $t \to \infty$.
\end{proof}

Thus $T_1$ and $T_2$ have natural density $0$.

%****************************************************************************
\section{Proof of Theorem~\ref{main}}

By Proposition~\ref{integerpoints}, $E'(\Z[S^{-1}])$ differs from 
$\{\,\pm \ell_i P: i \ge 1 \,\}$
by at most a finite set.
Since $y(\pm \ell_i P)$ is within $1/10$ of $\pm i$, 
any bounded subset of $\R^2$ contains at most finitely
many points of $E'(\Z[S^{-1}])$.
Part~(1) of Theorem~\ref{main} follows.

We next construct a diophantine model $A$ of the positive integers 
over $\Z[S^{-1}]$.
The set of nonzero elements of $\Z[S^{-1}]$ 
is diophantine (see Theorem 4.2~ of~\cite{shlapentokh1994}), 
and we can represent elements of $\Q$ as fractions of
elements of $\Z[S^{-1}]$ with nonzero denominator.
Therefore equations over $\Q$ can be rewritten as 
systems of equations over $\Z[S^{-1}]$,
and there is no harm in using them in our diophantine definitions.
In particular, we may use the predicate $x \ge y$,
since it can be encoded as 
$(\exists z_1,z_2,z_3,z_4 \in \Q)(x=y+z_1^2+z_2^2+z_3^2+z_4^2)$.

For $i \in \Z_{>0}$, define $y_i:=y(\ell_i P)$.
Let $A = \{y_1,y_2,\dots\}$.
Then $A$ is diophantine over $\Z[S^{-1}]$,
because it consists of the nonnegative elements
of the set of $y$-coordinates of $E'(\Z[S^{-1}])$ minus a finite set.
We have a bijection $\Z_{>0} \to A$ taking $i$ to $y_i$.

It remains to show that the graphs of addition and multiplication
on $\Z_{>0}$ correspond to diophantine subsets of $A^3$.
We know $|y_i-i| \le 1/(10i) \le 1/10$,
so the idea is that the addition on $\Z_{>0}$
should correspond to the operation of adding elements of $A$
and then rounding to the nearest element of $A$.
A similar idea will work for squaring,
and we will get multiplication from addition and squaring.

\begin{lemma}
\label{approximate}
Let $m,n,q \in \Z_{>0}$.
Then
\begin{enumerate}
\item $m+n=q$ if and only if $|y_m + y_n - y_q| \le 3/10$.
\item $m^2=n$ if and only if $|y_m^2 - y_n| \le 4/10$.
\end{enumerate}
\end{lemma}

\begin{proof}\hfill

(1) The quantity $y_m+y_n-y_q$ differs from the integer $m+n-q$ 
by at most $1/10+1/10+1/10$.

(2) The quantity $y_m^2 - y_n$ differs from the integer $m^2-n$
by at most 
\begin{align*}
	|y_m^2 - m^2| + |y_n - n| 
	&\le \left|\left(m+\frac{1}{10m}\right)^2 - m^2\right| + \frac{1}{10} \\
	&\le \frac{4}{10}.\qedhere
\end{align*}
\end{proof}

Lemma~\ref{approximate} shows that the two predicates
$m+n=q$ and $m^2=n$ on $\Z_{>0}$
correspond to diophantine predicates on $A$.
Building with these, we can show the same for $mn=q$,
since
	$$mn=q \quad\iff\quad (m+n)^2 = m^2 + n^2 + q + q.$$
This completes the proof of part~(2) of Theorem~\ref{main}.
Part~(3) follows from~(2) and Matijasevi{\v{c}}'s Theorem.

%****************************************************************************
\section*{Acknowledgements}

I thank Thanases Pheidas and the referee for helpful comments,
in particular for simplifying the 
diophantine definition of multiplication at the end.
I thank also Gunther Cornelissen and Graham Everest for suggesting
references for Section~\ref{S:denominators},
and Alexandra Shlapentokh for suggesting some improvements in the exposition.

\bibliographystyle{amsalpha}
\bibliography{bjorn}

\end{document}